\documentclass[a4paper,12pt]{article}
\newtheorem{definition}{{\bf Definition}}[section]
\newtheorem{theorem}[definition]{{\bf Theorem}}
\newtheorem{lemma}[definition]{{\bf Lemma}}
\newtheorem{proposition}[definition]{{\bf Proposition}}
\newtheorem{corollary}[definition]{{\bf Corollary}}
\newtheorem{example}[definition]{{\bf Example}}

\usepackage{amscd,amsmath,amssymb,graphicx}

\begin{document}

\title{Factors generated by $C^*$-finitely correlated states}
\author{Hiromichi Ohno \\
Graduate School of Information Sciences, Tohoku University,\\
Aoba-ku, Sendai 980-8579, Japan}
\date{}
\maketitle

\begin{abstract}
We present several
equivalent conditions for $C^*$-finitely correlated states 
defined on the UHF algebras to be
factor states and consider the types of factors generated by them.
Subfactors generated by generalized quantum Markov chains defined on
the gauge-invariant parts of the UHF algebras are also discussed.

\end{abstract}

\section{Introduction}

The notion of quantum Markov chains was introduced by Accardi
in [\ref{acc}].
As special cases,
the notion of quantum Markov states was defined by Accardi and
Frigerio in [\ref{accardi1}] and
that of $C^*$-finitely correlated states
was discussed by Fannes, Nachtergaele and Werner [\ref{fannes1}].
Further discussions
on quantum Markov states are found in [\ref{accardi2}], [\ref{go}]
and [\ref{ohno}] for example.

In [\ref{fidaleo}],
Fidaleo and Mukhamedov showed that the von Neumann algebras generated by
faithful translation-invariant 
quantum Markov states are factors
of type ${\rm II}_1$ or type ${\rm III}_\lambda$ with $\lambda \in (0,1]$.
In the present paper we discuss the von Neumann algebras generated
by $C^*$-finitely correlated states.
In the case where the states are Markov states,
it is known ([\ref{go}], [\ref{ohno}] for example)
that the states are unique KMS states, and
the exact form of local density matrices is also known.
Hence, we can see that the von Neumann algebras are factors, and
the types of factors can be determined in terms of the local density matrices.
But, in the case where the states are $C^*$-finitely correlated states,
we have to find a different method.

A $C^*$-finitely correlated stete is a state on the UHF algebra
$\bigotimes_{\mathbb Z} M_d$ defined by a triplet
$({\mathfrak C}, E, \rho)$, where
${\mathfrak C}$ is a finite dimensional $C^*$-algebra,
$E$ is a completely positive map from $M_d \otimes {\mathfrak C}$
to ${\mathfrak C}$ and $\rho$ is a state on ${\mathfrak C}$.

In section 2, we show that a $C^*$-finitely
correlated state is a factor state if and only if
it satisfies the strong mixing property.
To see this, we look at the eigenvectors of $E(I \otimes \, \cdot \, )$
with eigenvalues of modulus $1$.
In section 3, we show that the factors generated by
$C^*$-finitely correlated states are of type
${\rm I}_\infty$ or type ${\rm II}_1$ or type ${\rm II}_\infty$ or
type ${\rm III}_\lambda$ for some $\lambda \in (0,1]$.

The notion of generalized quantum Markov chains was
introduced in [\ref{ohno2}, \ref{ohno3}].
Generalized quantum Markov chains 
extend
translation-invariant quantum Markov chains to those on AF algebras.
In [\ref{ohno3}], we considered the case where the AF algebras
are gauge-invariant $C^*$-algbras and
we proved the extendability theorem for
any generalized quantum Markov chain,
that is, the generalized quantum Markov chain is restriction of
a quantum Markov chain on the UHF algebra.

In section 4, we show that a faithful generalized quantum Markov chain
defined on the gauge-invariant $C^*$-algebra is a factor state as long as
the extended state is factorial.
Moreover, we present some examples of subfactors
indeced from generalized quantum Markov chains.

\section{Equivalent condition for factor}

Let
${\mathfrak B}_i =M_d = M_d({\mathbb C})$,
the $d \times d$ complex matrix algebra, for $i\in {\mathbb Z}$
and
${\mathfrak B}$ be the infinite $C^*$-tensor product
$ \bigotimes_{i \in {\mathbb Z}} {\mathfrak B}_i$.
We denote
${\mathfrak B}_\Lambda =\bigotimes_{n \in \Lambda} {\mathfrak B}_n$ 
for arbitrary
subset $\Lambda \subset {\mathbb Z}$. The translation $\gamma$ is 
the right shift on ${\mathfrak B}$. We write $\phi_{[1,n]}$ 
for the localization 
$\phi|{\mathfrak B}_{[1,n]}$.
The following definition is from [\ref{fannes1}].

%%%%%%%%%%%%%%%%%%%5definition of C^* finitely corr states 

\begin{definition}{\rm
A state $\phi$ on ${\mathfrak B}$ is called a {\it $C^*$-finitely correlated
state} if there exist a finite dimensional $C^*$-algebra ${\mathfrak C}$,
a completely positive map $E:M_d \otimes {\mathfrak C}
\to {\mathfrak C}$ and a state $\rho$ on ${\mathfrak C}$ such that
\[
\rho(E(I_d \otimes C)) = \rho (C)
\]
for all $C \in {\mathfrak C}$ and
\[
\phi(A_1 \otimes \cdots \otimes A_n)
=\rho(E(A_1 \otimes E(A_2 \otimes \cdots \otimes E(A_n \otimes
I_{\mathfrak C}) \cdots )))
\]
for all $A_1 ,\ldots , A_n \in M_d$.}
\end{definition}

Let $\phi$ be a $C^*$-finitely correlated state generated by the triplet
$( {\mathfrak C}, E, \rho)$.
For any $n \in {\mathbb N}$, we define the completely positive map
$E^{(n)}$ from ${\mathfrak B}_{[1,n]} \otimes {\mathfrak C}$ to
${\mathfrak C}$ by
\[
E^{(n)}(A_1 \otimes \cdots \otimes A_n \otimes C)
= E(A_1 \otimes E(A_2 \otimes \cdots \otimes E(A_n \otimes C) \cdots ))
\]
for all $A_1, \ldots , A_n \in M_d$ and $C \in {\mathfrak C}$.
We will also need the linear space ${\mathfrak C}_0 \subset {\mathfrak C}$ 
which is the smallest subspace of ${\mathfrak C}$ containing $I$ and
invariant under $E(A \otimes \, \cdot \, )$ for all $A \in M_d$.
Since ${\mathfrak C}$ is finite dimensional, there exists an 
integer $N$ such that
\[
{\mathfrak C}_0 = \{ E^{(N)}(A_{[1,N]} \otimes I) \,|\, A\in 
{\mathfrak B}_{[1,N]} \}.
\]
Moreover, we assume that the triplet $( {\mathfrak C}, E, \rho)$ is minimal,
that is, ${\mathfrak C}_0$ generates ${\mathfrak C}$ in the sense of algebra.

Let $( {\cal H}, \pi, \xi)$ be the GNS representation of $\phi$. 
Then, we can extend $\phi$ to $\pi({\mathfrak B})''$.
In the following, we omit $\pi$, if there is no confusion.
We want to show the condition that $\pi({\mathfrak B})''$ is a factor.
To this end, we introduce two subspaces of ${\mathfrak C}_0$.
We define the subspaces $L(E)$ and $L_1(E)$ by
\[
L(E) = \{ C \in {\mathfrak C} \,\, |\,\, E_I(C) = 
\lambda C , \,\, \lambda \in {\mathbb T} \}
\]
and
\[
L_1(E) = \{ C \in {\mathfrak C} \,\, |\,\, E_I(C) = C \},
\]
where $E_I = E(I\otimes \, \cdot \,)$.
$L_1(E)$ is the eignespace of $E_1$ with eigenvalue $1$ and
$L(E)$ is the space generated by
eigenspaces with eigenvalues of modulus $1$.
From [\ref{fannes2}], $L(E)$ and $L_1(E)$ are algebras 
containd in the center
of ${\mathfrak C}$. Moreover, 
there exists an integer $M$ such that $\lambda^M = 1$
for any eigenvalue $\lambda$ of $E_I$ with
modulus $1$.

The following argument is in [\ref{fannes2}].
For any minimal projection $P$ of $L_1(E)$,
we consider the algebra ${\mathfrak C}_P = P{\mathfrak C}P$.
Obviously, ${\mathfrak C} = \bigoplus {\mathfrak C}_P$,
where the sum is taken over all minimal projections in $L_1(E)$.
Since $E$ is a completely positive map, we have
$E(M_d \otimes {\mathfrak C}_P) \subset {\mathfrak C}_P$.
Therefore, we can define the restriction 
$E_P : M_d \otimes {\mathfrak C}_P \to {\mathfrak C}_P$.
We can assume $\rho(P) \neq 0$.
Then, with $\rho_P = \rho(P)^{-1} \rho | {\mathfrak C}_P$,
we have a triplet $({\mathfrak C}_P, E_P, \rho_P)$ generating
a $C^*$-finitely correlated state $\phi_P$.
A direct expression of $\phi_P$ is
\begin{equation}
\phi_P(A_1 \otimes \cdots \otimes A_n)
=\rho(P)^{-1} \rho(E(A_1 \otimes \cdots \otimes
E(A_n \otimes P) \cdots )) \label{L1Estate}
\end{equation}
for all $A_1, \ldots , A_n \in M_d$.
Then, we have the decomposition
\[
\phi = \sum \rho(P) \phi_P,
\]
where the sum is taken over all minimal projections in $L_1(E)$.

Let $\Pi$ denote the set of minimal projections in $L(E)$.
Then,
$E_I | \Pi$ defines a bijective map
from $\Pi$ to $\Pi$.
For any projection $Q$ in $\Pi$,
we have $E_I^M (Q) =Q$.
Hence, $Q$ is in $L_1(E_I^{(M)})$, 
where $E_I^{(M)} = E^{(M)}(I^{\otimes M} \otimes \, \cdot \,)$, and
we have a $C^*$-finitely correlated state $\phi_Q$
on a regrouped chain generated by the triplet
$({\mathfrak C}_Q, E_Q^{(M)}, \rho_Q)$,
where ${\mathfrak C}_Q$ and $\rho_Q$ are defined as above and
$E_Q^{(M)}$ is the completely positive map from 
$\bigotimes^M M_d \otimes {\mathfrak C}_Q$ to ${\mathfrak C}_Q$ defined by
\[
E_Q^{(M)} (A_1\otimes A_2 \otimes \cdots \otimes A_M \otimes C_Q)
=E(A_1 \otimes E(A_2 \otimes \cdots \otimes E(A_M \otimes C_Q) \cdots))
\]
for any $A_1, \ldots ,A_M \in M_d$ and $C_Q \in {\mathfrak C}_Q$.
A direct expression of $\phi_Q$ is
\begin{eqnarray}
\phi_Q (A_1 \otimes \cdots \otimes A_n)
=\rho(Q)^{-1} \rho(E^{(M)}(A_1 \otimes \cdots \otimes E^{(M)}
(A_n \otimes Q)
\cdots )) \label{LEstate}
\end{eqnarray}
for all $A_1, \ldots , A_n \in \bigotimes_{i=1}^M M_d$.
Then, we have the decomposition
\[
\phi = \sum_{Q \in \Pi} \rho(Q) \phi_Q.
\]
Moreover, $\phi_Q$ is strongly clustering for $\gamma^M$, that is,
\[
\lim_{n\to \infty} \phi_Q(A \gamma^{nM}(B)) =\phi(A)\phi(B)
\]
for all $A,B \in {\mathfrak B}$.
Indeed, we consider the Jordan decomposition of 
$(E^{(M)}_Q)_I = E^{(M)}_Q( I^{\otimes M} \otimes \, \cdot \, )$,
i.e., 
\[
(E^{(M)}_Q)_I = \sum_{\lambda} (\lambda P_\lambda +R_\lambda ),
\]
where the sum is taken over all eigenvalues, 
$P_\lambda P_{\lambda'} = \delta_{\lambda \lambda'} P_\lambda$ and
$R_\lambda$ is nilpotent with $P_\lambda R_{\lambda'}
=R_{\lambda'} P_\lambda = \delta_{\lambda \lambda'} R_\lambda$.
Since $\| (E^{(M)}_Q)_I \| \le 1$ and $(E^{(M)}_Q)_I$ has
trivial peripheral spectrum ([\ref{fannes2}]), i.e., the only eigenvector of
$(E_Q^{(M)})_I$ with eigenvalue of modulus $1$ is $Q$, $R_1 = 0$ and
$P_\lambda = R_\lambda = 0$ for
$\lambda$ with $|\lambda| \geq 1$ and $\lambda \neq 1$.
Hence, for any $\varepsilon >0$, there exists a number $m \in {\mathbb N}$
such that
$\| P_1 - (E^{(M)}_Q)_I^m \| < \varepsilon$.
Furthermore, for any $A \in {\mathfrak B}_{[1,nM]}$, we obtain
\begin{eqnarray*}
\phi_Q(A) = \rho_Q (E^{(nM)}(A \otimes Q))
=\lim_{l\to\infty}  \rho_Q(
 (E^{(M)}_Q)_I^l(E^{(nM)}(A \otimes Q))).
\end{eqnarray*}
Therefore, we have
\[
\lim_{l\to\infty} (E^{(M)}_Q)_I^l(E^{(nM)}(A \otimes Q)) =\phi_Q(A) Q.
\]
This implies that $\phi_Q$ is strongly clustering for $\gamma^M$.
In particular, if $\Pi = \{ I \}$, we obtain
\begin{eqnarray}
\lim_{l \to \infty}( E_I^l (E^{(n)} (A \otimes I)) = \phi (A) I \label{eq2.1}
\end{eqnarray}
for all $A \in {\mathfrak B}_{[1,n]}$.

For each $Q \in \Pi$, we set the projection $\bar Q \in L(E)$ by
\[
\bar Q = \sum \{ R \in \Pi \,\, | \,\, \phi_Q = \phi_R \}
\]
and the set $\bar \Pi$ by
\[
 \bar \Pi = \{ \bar Q \,\, |\,\, Q \in \Pi \}.
\]

%%%%%%%%%%%%%%    lemma   L(E) \cap C_0

\begin{lemma}\label{L(E)}
With the above notation, we have
\[
L(E) \cap {\mathfrak C}_0 = {\rm span} \bar \Pi .
\]
\end{lemma}
\begin{proof}
For any $T \in L(E) \cap {\mathfrak C}_0$,
there exists an element $B \in {\mathfrak B}_{[1,nM]}$ such that
$E^{(nM)}(B \otimes I) = T$. From the above argument,
we have
\begin{eqnarray}
T &=& E^{(nM)}(B \otimes I) = 
\lim_{ l \to \infty} E_I^{lM}(E^{(nM)}(B \otimes I) ) \nonumber \\
&=& \lim_{ l \to \infty}
\sum_{Q \in \Pi} (E_Q^{(M)})_I^l (E_Q^{(nM)}(B \otimes Q)) \nonumber \\
&=&
\sum_{Q \in \Pi} \phi_Q (B) Q = \sum_{\bar Q \in \bar \Pi} 
\phi_{Q}(B) \bar Q. \label{eq2.2}
\end{eqnarray}
This implies 
$L(E) \cap {\mathfrak C}_0 \subset {\rm span} \bar \Pi.$ 

To prove the converse, we show that $\bar Q \in {\mathfrak C}_0$ for
any $Q \in \Pi$.
For each $P,Q \in \Pi$, $\bar P \neq \bar Q$ implies
$\phi_P \neq \phi_Q$.
Since $\phi_P$ and $\phi_Q$ are $\gamma^M$-ergodic,
$\phi_P \neq \phi_Q$ implies that $\phi_P$ and $\phi_Q$ are
mutually disjoint ([\ref{BR1}, 4.3.19]).
Hence, for any $\varepsilon > 0$, 
there exists
 an element $A \in {\mathfrak B}_{[-nM+1,nM]}$ such that
$| \phi_P (A) - 1 | < \varepsilon$
 and $| \phi_Q (A)| < \varepsilon$ for any $Q \in \Pi$ with 
$\bar P \neq \bar Q$.
Since $\phi_Q$ is $\gamma^M$-invariant, we can assume that
$A \in {\mathfrak B}_{[1,nM]}$ for some $n \in {\mathbb N}$. 
Moreover, from (\ref{eq2.2}), 
there exists a number $l \in {\mathbb N}$ such that
\[
\|  E_I^{lM} (E^{(nM)} (A\otimes I)) - 
\sum_{\bar Q \in \bar \Pi} \phi_{Q}(A) \bar Q \| < \varepsilon.
\]
Therefore we have
\begin{eqnarray*}
 && \| \bar P -  E_I^{lM} (E^{(nM)} (A\otimes I)) \| \\
&\le& 
\| \bar P - \sum_{\bar Q \in \bar \Pi} \phi_{Q}(A) \bar Q \|
+ \| \sum_{\bar Q \in \bar \Pi} \phi_{Q}(A) \bar Q
-E_I^{lM} (E^{(nM)} (A\otimes I)) \| \\
&<& 2 \varepsilon
\end{eqnarray*}
Since ${\mathfrak C}_0$ is closed and $E_I^{lM} (E^{(nM)} (A\otimes I))$
is in ${\mathfrak C}_0$, 
we have $\bar P \in {\mathfrak C}_0$.
\end{proof}

Now, we have the next theorem.

%%%%%%%%%%%%%%%%%%%%5     theorem    factor equivalent condition

\begin{theorem}\label{eqivalence}
For any $C^*$-finitely correlated state $\phi$ generated by
the triplet $({\mathfrak C}, E, \rho)$, the following
conditions are equivalent.\\
{\rm (i)} $\pi({\mathfrak B})''$ is a factor. \\
{\rm (ii)} $\phi$ is strongly clustering for $\gamma$. \\
{\rm (iii)}
$L(E) \cap {\mathfrak C}_0 = {\mathbb C} I$. \\
{\rm (iv)}
$\bar \Pi= \{ I \}$, that is, $\phi_Q$'s in {\rm (\ref{LEstate})} with 
projections $Q \in \Pi$ are same.
\end{theorem}
\begin{proof}
(iii) $\Leftrightarrow$ (iv) follows from Lemma \ref{L(E)}.

(iii) $\Rightarrow$ (ii). 
Since $L(E) \cap {\mathfrak C}_0 ={\mathbb C} I$ implies 
$\phi = \phi_Q$ for any $Q \in \Pi$,
$\phi$ is strongly clustering for $\gamma^M$.
Moreover, $\phi$ is $\gamma$-invariant.
Therefore, we have
\begin{eqnarray*}
\lim_{n \to \infty} \phi(A\gamma^{nM +l} (B))
=\lim_{n \to \infty} \phi(A \gamma^{nM}(\gamma^l(B))) =
\phi(A) \phi(\gamma^l(B)) = \phi (A) \phi(B)
\end{eqnarray*}
for any $A, B \in {\mathfrak B}$ and $0 \le l \le k-1$.
Hence, $\phi$ is strongly clustering for $\gamma$.

(i) $\Rightarrow$ (iv). 
For any $P, Q \in \Pi$, $\bar P \neq \bar Q$ implies
$\phi_P$ and $\phi_Q$ are disjoint.
This contradicts $Z(\pi({\mathfrak B})'') = {\mathbb C} I$.
Hence, we obtain $\bar \Pi = \{ I \}$.

%Since $E_I^k(Q) =Q$, for any element $Q \in L(E) \cap {\mathfrak C}_0$ and
%$B_n, B_n' \in {\mathfrak B}_{[1,n]}$ and $qk > pk > n$, we have
%\begin{eqnarray*}
%\langle B_n \xi, \gamma^{pk}(\eta(Q)) B_n' \xi \rangle
%&=& \phi(B_n^* \gamma^{pk}(\eta(Q)) B_n') 
%=\phi(B_n^* \gamma^{qk}(\eta(Q)) B_n') \\
%&=& \langle B_n \xi, \gamma^{qk}(\eta(Q)) B_n' \xi \rangle .
%\end{eqnarray*}
%Combining $Z(\pi({\mathfrak B})'') = {\mathbb C} I$ and Lemma \ref{atinfty},
%we get $\lim_{r\to \infty} \gamma^{rk}(\eta(Q)) = \alpha I$ for some
%$\alpha \in {\mathbb C}$.
%Moreover, we have
%\begin{eqnarray*}
%\lim_{r \to \infty} \langle \xi, \gamma^{rk}(\eta(Q)) \xi \rangle
%&=&\lim_{r \to \infty} \phi ( \gamma^{rk}(\eta(Q)) )
%=\phi(\eta(Q)) \\
%&=& \rho (E^{(N)} (\eta(Q) \otimes I))
%= \rho (Q).
%\end{eqnarray*}
%Hence, we obtain $\alpha = \rho(Q)$.
%Therefore, we have
%\begin{eqnarray*}
%\langle B_n \xi, \gamma^{pk}(\eta(Q)) B_n' \xi \rangle
%= \lim_{r \to \infty} \langle B_n \xi, \gamma^{rk}(\eta(Q)) B_n' \xi \rangle 
%= \rho(Q) \phi (B_n^* B_n').
%\end{eqnarray*}
%On the other hand, we obtain
%\begin{eqnarray*}
%&& \langle B_n \xi, \gamma^{pk}(\eta(Q)) B_n' \xi \rangle
%= \phi (B_n^* \gamma^{pk}(\eta(Q)) B_n')  \\
%&=& \rho (E^{(pk)}(B_n^* B_n' \otimes I^{\otimes pk-n} \otimes
%E^{(N)}(\eta(Q) \otimes I))) \\
%&=& \rho (E^{(pk)}(B_n^* B_n' \otimes I^{\otimes pk-n} \otimes Q))
%= \rho(Q) \phi_Q (B_n^* B_n').
%\end{eqnarray*}
%This means $\phi_Q = \phi$. Hence, we have
%$\bar \Pi =\{ I \}$. 

(ii) $\Rightarrow$ (iii).
We assume that $\phi$ is strongly clustering for $\gamma$.
Then, $\phi$ is strongly clustering for $\gamma^M$ and
hence $\gamma^M$-ergodic.
Since $\phi_Q$ is $\gamma^M$-ergodic
for any $Q \in \Pi$, we have $\bar\Pi = \{ I \}$.

(ii) $\Rightarrow$ (i).
Since $Z(\pi({\mathfrak B})'') = \bigcap_{n \in {\mathbb N}}
\pi({\mathfrak B}_{(-\infty, -n] \cup [n, \infty )} )''$ 
(see e.g. [\ref{BR1}, 2.6.10] ),
for any $X \in Z(\pi({\mathfrak B})'')$ with $\| X\| =1$,
there exists a sequence
$\{ X_n  \}$ with $X_n \in {\mathfrak B}_{[-l(n),n] \cup [n,l(n)]}$, 
$\| X_n \| \le 1$ and 
$\lim_{n\to\infty} X_n = X$ in the weak operator topology.
We can write 
\[
X_n = \sum Y_i^{(n)} \gamma^{n-1} ( Z_i^{(n)})
\]
for some $Y_i^{(n)} \in {\mathfrak B}_{[-l(n), -n]}$ and
$Z_i^{(n)} \in {\mathfrak B}_{[1, l(n)-n+1]}$.
For any element $A \in {\mathfrak B}_{[1,p]}$, $p \in {\mathbb N}$, 
there exists an
element $A' \in {\mathfrak B}_{[1,N]}$ such that
\[
E^{(p)}(A \otimes I) = E^{(N)} (A' \otimes I).
\]
We write $A' = \theta(A)$.
For any element $B_m, B_m' \in {\mathfrak B}_{[1,m]}$ with $m<n$,
we have
\begin{eqnarray*}
&& \langle B_m \xi ,  (I^{\otimes n} \otimes A ) B_m' \xi \rangle
= \phi ( B_m^* (I^{\otimes n} \otimes A ) B_m' ) \\
&=& \rho ( E^{(n)} (B_m^* B_m' \otimes I^{\otimes n-m} \otimes 
E^{(p)}( A \otimes I))) \\
&=&\rho ( E^{(n)} (B_m^* B_m' \otimes I^{\otimes n-m} \otimes 
E^{(N)}( \theta(A) \otimes I))) \\
&=&  \langle B_m \xi ,  (I^{\otimes n} \otimes \theta(A) ) B_m' \xi \rangle.
\end{eqnarray*}
Therefore, $X_n' = \sum Y_i^{(n)} \gamma^{n-1} ( \theta(Z_i^{(n)}))$
converges to $X$ in the weak operator topology.
Moreover, since $\theta(Z_i^{(n)})  \in  {\mathfrak B}_{[1,N]}$,
we can write
\[
X_n' = \sum_{i=1}^{d^{2N}} S_i^{(n)} \gamma^n (T_i)
\]
for some $S_i^{(n)} \in {\mathfrak B}_{[-l(n), -n]}$ and
a system of matrix units $\{ T_i\}$ of ${\mathfrak B}_{[1,N]}$.
Since $X_n'$ converges to $X$ in the weak operator topology,
there exists some constant $C>0$ such that $\| X_n' \| \le C$ for any
$n \in {\mathbb N}$.
Then, we have $\| S_i^{(n)} \| \le C$.

From the proof of (\ref{eq2.1}), for $\varepsilon >0$
there exists $L \in {\mathbb N}$ such that
\[
\| E_I^L (E^{(p)}(A\otimes I)) - \phi(A) I \| < \varepsilon \| A\|
\]
for any $A \in {\mathfrak B}_{[1,p]}$ and $p \in {\mathbb N}$.
Using this uniform convergence, 
for any $B_m, B_m' \in {\mathfrak B}_{[1,m]}$ we have
\begin{eqnarray*}
&& 
\langle B_m \xi , X B_m' \xi \rangle 
=
\lim_{n\to\infty}
\langle B_m \xi , X_n' B_m' \xi \rangle \\
&=&\lim_{n\to\infty} \sum_{i=1}^{d^{2N}} 
\langle B_m \xi , S_i^{(n)} \gamma^n (T_i) B_m' \xi \rangle 
=\lim_{n\to \infty} \sum_{i=1}^{d^{2N}}  
\phi(B_m^*  S_i^{(n)} \gamma^n (T_i)  B_m') \\
&=& \lim_{n\to \infty} \sum_{i=1}^{d^{2N}}  
\rho( E^{l(n)-n + 1}(S_i^{(n)} \otimes E_I^n ( E^{(m)}(B_m^*B_m' \otimes
E_I^{n-m}(E^{(N)}(T_i \otimes I)))))) \\
&=& \lim_{n\to\infty} \sum_{i=1}^{d^{2N}}  \phi(S_i^{(n)})
\phi(B_m^*B_m') \phi(T_i) 
= \phi(B_m^*B_m') \lim_{n\to\infty} \sum_{i=1}^{d^{2N}} 
\phi( S_i^{(n)} \gamma^n (T_i)) \\
&=& \phi(B_m^*B_m') \lim_{n\to\infty} \phi(X_n') 
=\langle B_m \xi, \phi(X) B_m' \xi \rangle.
\end{eqnarray*}
Therefore, we obtain $X = \phi(X) I$.
\end{proof}

By the theorem, 
for any $ P , Q \in \Pi$ such that $\phi_P \neq \phi_Q$,
$\phi_P$ and $\phi_Q$ are disjoint and factor states.
Therefore, for any $P \in \Pi$,
there exists a minimal projection $T$ in $Z(\pi({\mathfrak B})'')$,
such that
\[
 \phi_P(B)=
\langle \xi ,T\xi \rangle^{-1} \langle \xi, B T\xi \rangle 
\]
for any $B \in \pi({\mathfrak B})''$.
In fact, $T$ is the support projection of $\phi_P$.
We define a bijective map $\eta$ from $\bar\Pi$ to 
a set of minimal projections in $Z(\pi({\mathfrak B})'')$ by
\[
\eta(\bar P) = T.
\]
Now, we have the next corollary.

%%%%%%%%%%%%%%%%%%%%    cor center

\begin{corollary}
We obtain
\[
Z(\pi({\mathfrak B})'') = {\rm span}\{ \eta(\bar P) \,\,
| \,\, \bar P \in \bar\Pi \}.
\]
In particular, the dimension of the center $Z(\pi({\mathfrak B})'')$ is
finite and not greater than the dimension of the center of ${\mathfrak C}$.
\end{corollary}

%%%%%%%%%%%%%%%%%%%%%%%%%%%%%%%% sect    examples

\section{Types of factors generated by $C^*$-finitely correlated states}

In this section, we examine the types of factors generated by
strongly clustering $C^*$-finitely correlated states.
In the following, we assume that
$\phi$ is a $C^*$-finitely correlated state generated by a
triplet $({\mathfrak C}, E, \rho)$ and
it is strongly clustering.

Since $\phi$ is $\gamma$-invariant, we can extend $\gamma$ to
$\pi({\mathfrak B})''$.
Let $P$ be the support projection of $\phi$.
Then, $\gamma (P) = P$. Indeed, 
$\phi(\gamma (P)) = \phi(P)$ implies
$\gamma (P) \ge P$. Similary, we have
$\gamma^{-1} (P) \ge P$.
This means $\gamma (P) = P$.
Therefore, we can define the automorphism $\gamma|P{\mathfrak B}P$.
Here, the normal extension of $\phi$ to $\pi({\mathfrak B})''$ is 
denoted by the same $\phi$ and $\pi({\mathfrak B})$ is identified with
${\mathfrak B}$.

Let $S(\pi({\mathfrak B})'')$ be the Connes invariant.
The next proposition is in [\ref{fidaleo}].
The proof is given for convenience.

%%%%%%%%%%%%%%%%%%%%55    lemma  about S set

\begin{proposition}
Let $\phi^P = \phi | P {\mathfrak B} P$. Then, we have
\[
S(\pi({\mathfrak B})'')\backslash \{ 0\} 
= {\rm Sp}(\Delta_{\phi^P}) \backslash \{0\},
\]
where $\Delta_{\phi^P}$ is a modular operator of $\phi^P$.
\end{proposition}
\begin{proof}
Since $\pi({\mathfrak B})''$ is a factor,
we know that $S(\pi({\mathfrak B})'') =S(P\pi({\mathfrak B})''P)$.
$P{\mathfrak B}P$ is asymptotically abelian with respect to $\gamma$
and $\phi^P$ is strongly clustering for $\gamma$.
Therefore, if a state $\omega$ on $P{\mathfrak B}P$ is quasi-containd
in $\phi^P$, then we have ${\rm Sp}(\Delta_{\phi^P})\subset
{\rm Sp}(\Delta_\omega)$ ([\ref{stormer}]).
In particular, for a projection $Q\in \pi({\mathfrak B})''$ with 
$0 \neq Q \le P$,
we have ${\rm Sp}(\Delta_{\phi^P}) \subset {\rm Sp}(\Delta_{\phi^Q})$,
where $\phi^Q = \phi^P(Q)^{-1} \phi^P(Q \,\, \cdot \,\, )$.
Moreover, $\phi^P$ is faithful on $P\pi({\mathfrak B})''P$
and $P{\mathfrak B}P$ is weakly dense in $P\pi({\mathfrak B})''P$.
Hence, we have 
\[
S(P\pi({\mathfrak B})''P)\backslash\{0\} =
{\rm Sp}(\Delta_{\phi^P}) \backslash \{0\}.
\]
\end{proof}

In the following, we examine the type of $\pi({\mathfrak B})''$.
In the case where ${\rm Sp}(\Delta_{\phi^P}) \neq \{	1 \}$,
since $\phi^P$ is faithful, ${\rm Sp}(\Delta_{\phi^P})$ contains
a number which is neither $0$ nor $1$.
Therefore, $S(\pi({\mathfrak B})'') \neq \{0,1\}$.
Hence, $\pi({\mathfrak B})''$ is a ${\rm III}_\lambda$ factor for some
$\lambda \in (0,1]$.

If ${\rm Sp}(\Delta_{\phi^P}) = \{ 1\}$,
then $\phi^P$ is a tracial state on $P\pi({\mathfrak B})''P$.
Hence, $P$ is a finite projection.
Therefore, $\pi({\mathfrak B})''$ is not a ${\rm III}$ factor.
If $\phi$ is faithful,
then $\pi({\mathfrak B})''$ is a ${\rm II}_1$ factor.
If $\phi$ is pure, then $\pi({\mathfrak B})''$ is a ${\rm I}_\infty$ factor.
From [\ref{fannes2}], $\phi$ is pure if and only if
$\phi$ is strongly clustering and the mean entropy of $\phi$ is zero.

%%%%%%%%%%%%%%%%%%%%%%%%%%%%%%5    prop II_\infty

\begin{proposition}
If ${\rm Sp}(\Delta_{\phi^P}) = \{ 1\}$ and $\phi$ is neither
locally faithful nor pure,
then $\pi({\mathfrak B})''$ is a ${\rm II}_\infty$ factor.
\end{proposition}
\begin{proof}
From the assumption, $\pi({\mathfrak B})''$ is either 
a ${\rm II}_1$ factor or a ${\rm II}_\infty$ factor.
Now, we assume that $\pi({\mathfrak B})''$ is a ${\rm II}_1$ factor.
Then, there is a faithful tracial state $\tau$ on $\pi({\mathfrak B})''$.
Since $\phi$ is not locally faithful, there exist a number 
$M \in {\mathbb N}$ and a projection $Q \in {\mathfrak B}_{[1,M]}$
such that $\phi(Q)=0$. 
We obtain $P \le I - Q$. Since $\gamma(P) = P$, we have
$P \le \gamma^{lM}(I - Q)$ for any $l \in {\mathbb N}$.
Let ${\rm rank}(Q) = q$ in ${\mathfrak B}_{[1,M]}$.
Then, we have
\begin{eqnarray*}
\tau(P) \le \tau ( \prod_{l=1}^n \gamma^{lM} (I - Q)) = 
(1- {q \over d^{M}})^n
\end{eqnarray*}
for any $n \in {\mathbb N}$.
This implies $\tau(P) =0$. Hence, 
$\pi({\mathfrak B})''$ is a ${\rm II}_\infty$ factor.
\end{proof}

In the rest of this section, we present examples of
${\rm III}_\lambda$ factors for $\lambda \in (0,1]$
which are generated by
translation-invariant quantum Markov states.

%%%%%%%%%%%%%%%%%%%%%%%%%%%%%%55     def of QMS

\begin{definition}{\rm [\ref{accardi1}]
A state $\phi$ on ${\mathfrak B}$ is said to be a 
{\it quantum Markov state},
if there exists a conditional expectation
$E_n$ from ${\mathfrak B}_{[1,n+1]}$ to ${\mathfrak B}_{[1,n]}$
such that ${\mathfrak B}_{[1,n-1]} \subset {\rm ran}(E_n)$ and
\[
\phi \circ E_n = \phi_{[1,n+1]}
\]
for each $n \in {\mathbb N}$.
}
\end{definition}

Although the above definition is a bit different from the original one 
of Accardi and Frigerio in [\ref{accardi1}], it is known that both
definitions are equivalent 
([\ref{go}]).

In the case where the quantum Markov state $\phi$ is translation-invariant,
we can assume that $E_n = {\rm id}_{{\mathfrak B}_{[1,n-1]}} \otimes E$
for some conditional expectation $E$ from $M_d \otimes M_d$ into $M_d$ 
([\ref{ohno}]).
Therefore, translation-invariant quantum Markov states are
$C^*$-finitely correlated states.

In the following, we assume that
$\phi$ is a locally faithful translation-invariant quantum Markov state
generated by $(E, \rho)$ with $\rho = \phi | {\mathfrak B}_1$ 
and that $\phi$ is not a tracial state.
Let ${\mathfrak D} = {\rm ran}(E)$. Since ${\mathfrak D}$ is a
finite dimensional $C^*$-algebra, we can write
\[
{\mathfrak D} = \bigoplus_{i=1}^p M_{d_i}.
\]
Let $m_i$ be the multiplicity of $M_{d_i}$ as a
$C^*$-subalgebra of $M_d$, and
we define
\[
\bar{\mathfrak D} = \bigoplus_{i=1}^p M_{m_i},
\]
${\mathfrak E}_n = \bar{\mathfrak D} \otimes {\mathfrak B}_{[1,n-1]} 
\otimes {\mathfrak D}$ and
${\mathfrak E}_n^{xy} = M_{m_x} \otimes {\mathfrak B}_{[1,n-1]} \otimes
M_{d_y}$ for $1 \le x,y \le p$.
From [\ref{accardi2}], there exist positive operators
$T_{ij} \in M_{m_i} \otimes M_{d_j}$ for any $1 \le i,j \le p$
such that the density matrix of 
$\phi | {\mathfrak E}_n$
is written by
\begin{eqnarray}
D_n = \bigoplus_{i_1, \ldots , i_n} 
\rho(I_{m_{i_1}})
T_{i_1i_2} \otimes 
T_{i_2i_3} \otimes \cdots \otimes
T_{i_{n-1}i_n}.\label{density1}
\end{eqnarray}
Since $T_{ij}$ is positive, we can choose a system of 
matrix units $\{e_{kl}^{(ij)} \}$ for $M_{m_i}\otimes M_{d_j}$
and write
\[
T_{ij}= {\rm diag} ( e^{t_1^{(ij)}}, e^{t_2^{(ij)}},
 \ldots , e^{t_{m_i d_j}^{(ij)}}).
\]
To calculate $S(\pi({\mathfrak B})'')$,
we consider ${\rm sp}(\Delta_\phi)$.
Since $\phi$ is faithful, we obtain
\[
{\rm sp}(\Delta_\phi) \backslash \{ 0\} = \exp ({\rm sp}(\sigma^\phi)),
\]
where $\sigma^\phi$ is the modular automorphism group of $\phi$
and ${\rm sp}(\sigma^\phi)$ is the Arveson spectrum of $\sigma^\phi$.
Since ${\mathfrak B}$ is weakly dense in $\pi({\mathfrak B})''$, we have
\begin{eqnarray*}
{\rm sp}(\sigma^\phi)
&=& \overline{\bigcup_{B\in {\mathfrak B}} {\rm sp}_{\sigma^\phi}(B) }
=\overline{\bigcup_{n=1}^\infty \bigcup_{B \in
{\mathfrak E}_n}  {\rm sp}_{\sigma^\phi}(B) } \\
&=& \overline{\bigcup_{n=1}^\infty \bigcup_{x,y=1}^p \bigcup_{B \in
{\mathfrak E}_n^{xy}}
  {\rm sp}_{\sigma^\phi}(B) }.
\end{eqnarray*}
From [\ref{accardi1}], we know that
\[
\sigma_t^\phi | {\mathfrak E}_n = {\rm Ad} D_n^{it}.
\]
Therefore, ${\mathfrak E}_n^{xy}$
is invariant under $\sigma^\phi$ and we have
\[
\bigcup_{B \in
{\mathfrak E}_n^{xy}}
  {\rm sp}_{\sigma^\phi}(B)
  = {\rm sp}(\sigma^\phi | 
{\mathfrak E}_n^{xy}).
\]

%%%%%%%%%%%%%%%%%%%%%%%%%%%%%%%%%5555    lemma

\begin{lemma}\label{spec1}
Let $\psi$ be a state on $M_k$ with the density matrix
$D = {\rm diag} (e^{t_1} ,\ldots ,e^{t_k})$.
Then the Arveson spectrum of $\sigma^\psi$ is written as
\[
{\rm sp}(\sigma^\psi) = \{ t_i - t_j \,\, | \,\, 1\le i,j \le k \}.
\]
\end{lemma}
\begin{proof}
This is obvious from the fact that
\[
\sigma^\psi_t = {\rm Ad}(D^{it}).
\]
\end{proof}

Since the density matrix of $\phi | {\mathfrak E}_n$ is written as in
(\ref{density1}), the density matrix of 
$\phi | {\mathfrak E}_n^{xy}$ is
written as
\[
\bigoplus_{i_2, \ldots , i_{n-1}} \rho(I_{m_x})
T_{x i_2} \otimes T_{i_2 i_3}
\otimes \cdots \otimes T_{i_{n-2}i_{n-1}} \otimes T_{i_{n-1}y}.
\]
Therefore, we have 
\begin{eqnarray}
&& {\rm sp}(\sigma^\phi | {\mathfrak E}_n^{xy}) \nonumber \\
&=&
\{ t_{q_1}^{(x i_2)} + \sum_{k=2}^{n-2} t_{q_k}^{(i_ki_{k+1})}
+ t_{q_{n-1}}^{(i_{n-1} y)}
-t_{r_1}^{(x j_2)} - \sum_{l=2}^{n-2} t_{r_l}^{(j_l j_{l+1})}
- t_{r_{n-1}}^{(j_{n-1} y)} \nonumber \\
&&  | \,\, {\rm all \,\, possible \,\, } \label{eq3.4}
i_k,j_l,q_k,r_l \}.
\end{eqnarray}

Since $\exp({{\rm sp}(\sigma^\phi)}) = 
S(\pi({\mathfrak B})'')\backslash \{ 0 \}$, 
${\rm sp}(\sigma^\phi) $ is a group.
Hence, we obtain 
${\rm sp}(\sigma^\phi) = {\mathbb R}$ or else
there exists a number $\lambda \in (0,1)$ such that
\[
{\rm sp}(\sigma^\phi) = (\log \lambda) {\mathbb Z}.
\]
Let $G$ be a closed subgroup of ${\mathbb R}$ generated by
\[
\{ t_{j_1}^{(i_1i_2)} + t_{j_2}^{(i_2i_4)}
 - t_{j_3}^{(i_1i_3)} - t_{j_4}^{(i_3i_4)} \,\,|\,\,
{\rm all \,\, possible \,\,} i_k, j_l \}.
\]

%%%%%%%%%%%%%%%%%%%%%%%5   prop G

\begin{proposition}\label{type1}
We obtain
\[
G = {\rm sp}(\sigma^\phi) 
\]
\end{proposition}
\begin{proof}
By (\ref{eq3.4}), for any $i_k, j_l$, we obtain 
\[
t_{j_1}^{(i_1i_2)} + t_{j_2}^{(i_2i_4)}
 - t_{j_3}^{(i_1i_3)} - t_{j_4}^{(i_3i_4)} \in 
 {\rm sp}(\sigma^\phi | {\mathfrak C}_2^{i_1i_4}).
\]
Therefore, $G \subset {\rm sp}(\sigma^\phi)$.

We show the converse.
From definition, we obtain
$t_{j_1}^{(i_1 i_1)} - t_{j_4}^{(i_4i_4)} \in G$.
Then, for any
\[
t_{j_1}^{(xi_1)} + t_{j_2}^{(i_1i_2)} + t_{j_3}^{(i_2y)}
- t_{j_4}^{(xi_3)} - t_{j_5}^{(i_3i_4)}- t_{j_6}^{(i_4 y)}
\in {\rm sp}(\sigma^\phi | {\mathfrak C}_3^{xy}),
\]
by adding 
\begin{eqnarray*}
&& t_{k_1}^{(i_1 i_1)}
+t_{j_5}^{(i_3i_4)} -t_{k_2}^{(i_3i_1)} - t_{k_3}^{(i_1i_4)} \\
&=& (t_{k_1}^{(i_1 i_1)} - t_{k_4}^{(i_4i_4)})
+(t_{j_5}^{(i_3i_4)} +t_{k_4}^{(i_4i_4)}
-t_{k_2}^{(i_3i_1)} - t_{k_3}^{(i_1i_4)}) \in G,
\end{eqnarray*}
we have
\begin{eqnarray*}
&& (t_{j_1}^{(xi_1)} + t_{j_2}^{(i_1i_2)} + t_{j_3}^{(i_2y)}
- t_{j_4}^{(xi_3)} - t_{j_5}^{(i_3i_4)}- t_{j_6}^{(i_4 y)}) \\
&+&
(t_{k_1}^{(i_1 i_1)}
+t_{j_5}^{(i_3i_4)} -t_{k_2}^{(i_3i_1)} - t_{k_3}^{(i_1i_4)}) \\
&=& (t_{j_1}^{(xi_1)}+ t_{k_1}^{(i_1 i_1)}
-t_{j_4}^{(xi_3)} -t_{k_2}^{(i_3i_1)} )
+ (t_{j_2}^{(i_1i_2)} + t_{j_3}^{(i_2y)}
-t_{k_3}^{(i_1i_4)} - t_{j_6}^{(i_4 y)}) \in G.
\end{eqnarray*}
Hence, we get ${\rm sp}(\sigma^\phi | {\mathfrak C}_3^{xy}) 
\subset G$.
The idea of the above calculation is to split $(xi_1i_2y, xi_3 i_4 y)$
to $(x i_1 i_1, x i_3 i_1)$ and $(i_1 i_2 y, i_1 i_4 y)$.
The same can be applied to longer words.
For example, split $(x i_1i_2i_3y,xi_4i_5i_6y)$ to
$(x i_1 i_1, x i_4 i_1)$, $(i_1 i_2 i_1, i_1 i_5 i_1)$ and
$(i_1 i_3 y, i_1 i_6 y)$.
In this way, we obtain ${\rm sp}(\sigma^\phi | {\mathfrak C}_n^{xy})$
for all $1 \le x,y \le p$ and $n \in {\mathbb N}$,
so that ${\rm sp}(\sigma^\phi) \subset G$.
\end{proof}

Now, we define a number $\lambda \in {\mathbb R}$ to be 
$1$ if $G = {\mathbb R}$ or to be
$t$ if $G = (\log t){\mathbb Z}$.
Then, we have the next proposition. 

%%%%%%%%%%%%%%%%%%%%%%%%%%%%%%   prop main

\begin{proposition}\label{main}
With the above definition, if $\phi$ is not a tracial state,
$\pi({\mathfrak B})''$ is a type ${\rm III}_\lambda$ factor.
\end{proposition}

It was shown in [\ref{fidaleo}] that
$\pi({\mathfrak B})''$ is a type ${\rm III}_\lambda$ factor
for some $\lambda \in (0,1]$ as far as $\phi$ is not tracial.
But, the above proposition enables us to determine
the $\lambda$ from the
density matrices $T_{ij}$'s.

%%%%%%%%%%%%%%%%%%%%%%%%%%%%%%%%%%%%%%%     sect 4

\section{Subfactors generated by generalized quantum Markov chains}

In this section, we consider the types of factors generated by
generalized Markov chain.
The notion of generalized Markov chains was
introduced in [\ref{ohno3}] and they are defined on AF algebras.
Since it is not easy to
treat the general case of AF algebras in the one-dimensional lattice,
we restrict our consideration to gauge-invariant
$C^*$-algebras of UHF algebras.
Let $G$ be a unitary subgroup of ${\cal U}(M_d)$.
Then, for any $g\in G$, we can define an automorphism 
$\alpha_g = {\displaystyle\lim_{\longrightarrow}}
{\rm Ad} g^n$ on ${\mathfrak B}$.
We set
\[
{\mathfrak A}_{[-n,m]} = ({\mathfrak B}_{[-n,m]})^G
=\{ B \in {\mathfrak B}_{[-n,m]} \,\, | \,\, 
\alpha_g ( B )  = B, \,\, g \in G \}.
\]
We define
\[
{\mathfrak A} = \overline{\bigcup_{n=1}^\infty {\mathfrak A}_{[-n,n]}}.
\]
Then, ${\mathfrak A}$ is an AF algebra and $\gamma | {\mathfrak A}$ is 
an automorphism on ${\mathfrak A}$.

%%%%%%%%%%%%%%%%%%%%%%%%%%%   def generalized Markov state

\begin{definition}
{\rm
Let $\tilde\phi$ be a $C^*$-finitely correlated state
generated by the triplet $( M_d , E, \rho)$.
We assume that $E$ satisfies the $G$-covariant condition,
that is,
\begin{eqnarray}
E ( g^{\otimes 2} (A \otimes B) g^{\otimes 2 *} ) =
g E (A\otimes B) g^* \label{eq4.1}
\end{eqnarray}
for all $A,B \in M_d$.
Then, we have
\[
({\rm id}_{{\mathfrak B}_{[-m,n-1]}} \otimes E) ({\mathfrak A}_{[-m, n+1]})
\subset {\mathfrak A}_{[-m,n]}.
\]

We call the state $\phi = \tilde\phi |{\mathfrak A}$ a
{\it generalized Markov chain} on ${\mathfrak A}$
generated by $E$ and $\rho$.
}
\end{definition}

Since we have
\[
({\rm id}_{{\mathfrak B}_{[1,n-1]}} \otimes E) ({\mathfrak A}_{[1, n+1]})
\subset {\mathfrak A}_{[1,n]},
\]
by putting $E_n = {\rm id}_{{\mathfrak B}_{[1,n-1]}} 
\otimes E| {\mathfrak A}_{[1,n+1]}$, we obtain
\[
\phi(A) = \rho \circ E_1 \circ E_2 \circ \cdots \circ E_{n}(A\otimes I)
\]
for any $A \in {\mathfrak A}_{[1,n]}$.
The above formula justifies the terminology for $\phi$.
Moreover, the extendability theorem in [\ref{ohno3}] says that
any generalized Markov chain (in the sense of [\ref{ohno3}])
on ${\mathfrak A}$ can be written as the restriction of some
$C^*$-finitely correlated state on ${\mathfrak B}$ as above.

In the following, Let $E$, $\rho$, $\phi$ and $\tilde\phi$
be given as above and $({\cal H}, \pi , \xi)$ be the
GNS representation of ${\mathfrak B}$ associated with
$\tilde\phi$. Moreover, we assume that
$G$ is a discrete group and that $\tilde\phi$ 
is faithful and strongly clustering.

Now, let us know that $\pi({\mathfrak A})''$ is a factor.
To this end, we need to show that
for any $g \in G \backslash \{ e \}$ 
the automorphism $\alpha_g$ (extended to $\pi({\mathfrak B})''$)
is outer
(see [\ref{stratila}], 22.3 and 22.14).
To see this, we first give two lemmas.

%%%%%%%%%%%%%%%%%%%%%%%%%55   lemma  invariant under \alpha_g

\begin{lemma}
$\tilde\phi$ is invariant under $\alpha_g$ for any $g \in G$.
\end{lemma}
\begin{proof}
Since $E$ satisfies $G$-covariant condition, we have
\begin{eqnarray*}
&& \phi(gA_1g^*\otimes gA_2 g^* \otimes \cdots \otimes 
gA_{n-1} g^* \otimes gA_ng^*)\\
&=& \rho(E(gA_1 g^* \otimes E (gA_2g^* \otimes \cdots \otimes 
E(gA_{n-1}g^*\otimes E(gA_ng^* \otimes I)) \cdots))) \\
&=& \rho(E(gA_1 g^* \otimes E (gA_2g^* \otimes \cdots \otimes 
E(gA_{n-1}g^*\otimes gE(A_n \otimes I)g^*) \cdots))) \\
&\vdots& \\
&=& \rho(gE(A_1  \otimes E (A_2 \otimes \cdots \otimes 
E(A_{n-1}\otimes E(A_n \otimes I)) \cdots))g^*)
\end{eqnarray*}
for any $A_1, \ldots ,A_n \in M_d$. Therefore, we have
\begin{eqnarray*}
&& \phi(gA_1g^*\otimes gA_2 g^* \otimes \cdots gA_{n-1} g^* \otimes gA_ng^*)\\
&=&\phi( I^{\otimes m}\otimes gA_1g^*\otimes gA_2 g^* \otimes 
\cdots gA_{n-1} g^* \otimes gA_ng^*) \\
&=&\lim_{m \to \infty} \phi
( I^{\otimes m}\otimes gA_1g^*\otimes gA_2 g^* \otimes 
\cdots gA_{n-1} g^* \otimes gA_ng^*) \\
&=&\lim_{m \to \infty} \rho(gE_I^m(E(A_1  \otimes E (A_2 \otimes \cdots
E(A_{n-1}\otimes E(A_n \otimes I)) \cdots)))g^*)\\
&=& \phi(A_1 \otimes A_2 \otimes \cdots \otimes A_n),
\end{eqnarray*}
where the last equation follows from the fact that
\[
\lim_{m \to \infty} E_I^m(E^{(k)}(A \otimes I))
=\phi (A)
\]
for all $A \in {\mathfrak B}_{[1,k]}$.
Hence, $\phi$ is invariant under $\alpha_g$.
\end{proof}

Since $\tilde\phi$ is invariant under $\alpha_g$, we can extend
$\alpha_g$ to $\pi({\mathfrak B})''$.

%%%%%%%%%%%%%%%%%%%%%%%%%%    lemma  gamma(U) = \lambda U

\begin{lemma}\label{implement}
If $\alpha_g$ is an inner automorphism and 
$U\in \pi({\mathfrak B})''$ is any implementing unitary of $\alpha_g$,
we have 
\[
\gamma(U) = \lambda U
\]
for some $\lambda \in {\mathbb C}$ with $|\lambda | =1$.
\end{lemma}
\begin{proof}
This is obvious since $\pi({\mathfrak B})''$ is a factor
by Theorem \ref{eqivalence}
%From definition, it is immediate that $\gamma \circ \alpha_g = \alpha_g
%\circ \gamma$.
%For any $B \in {\mathfrak B}$, we have
%\begin{eqnarray*}
%\gamma(U) B \gamma(U)^* &=&  \gamma( U \gamma^{-1}(B) U^*)
%=\gamma \circ \alpha_g \circ \gamma^{-1} (B) \\
%&=& \alpha_g (B).
%\end{eqnarray*}
%Therefore, $\gamma(U)$ is also an implementing unitary of $\alpha_g$.
%This implies that
%\begin{eqnarray*}
%\gamma(U)^* U B U^* \gamma(U) = \gamma(U)^* \alpha_g(B) \gamma(U) = B
%\end{eqnarray*}
%for any $B \in {\mathfrak B}$.
%Hence, we have $\gamma(U)^* U$ is in the center $Z(\pi({\mathfrak B})'')$.
%Moreover, since $\gamma(U)^* U$ is a unitary operator, we have
%\[
%\gamma(U) = \lambda U
%\]
%for some $\lambda \in {\mathbb C}$ with $|\lambda | =1$.
\end{proof}

Now, we have the next proposition.

%%%%%%%%%%%%%%%%%%%%%%%%%%%%%5    prop   outer auto

\begin{proposition}
$\alpha_g$ is an outer automorphism for any $g\in G \backslash \{ e\}$.
\end{proposition}
\begin{proof}
We assume that $\alpha_g$ is an inner automorphism.
Then, there exists a unitary operator $U\in \pi({\mathfrak B})''$ such that
$\alpha_g = {\rm Ad}U$.
By Lemma \ref{implement}, there exists $\lambda \in {\mathbb C}$ with
$|\lambda |=1$ such that $\gamma(U) =\lambda U$.
For any $\varepsilon >0$, there exist
$U_m \in {\mathfrak B}_{[-m,m]}$ such that
\[
\| U_m \xi -U \xi \| < \varepsilon
\] 
and a sufficiently large number 
$L \in {\mathbb N}$ such that
\[
|\lambda^L-1|<\varepsilon.
\]
Since $\tilde\phi$ satisfies strongly clustering and $L$ is sufficiently 
large,  
for any $B_n \in {\mathfrak B}_{[1,n]}$ with $\| B_n \| \le 1$,
we have
\begin{eqnarray*}
&& | \langle B_n \xi , U \xi \rangle - \tilde\phi(B_n^*) \tilde\phi(U)| \\
&\le& |\langle B_n \xi , U \xi \rangle - 
\langle B_n \xi, \gamma^L(U) \xi \rangle |
+|\langle B_n \xi, \gamma^L(U) \xi \rangle -
\tilde\phi(B_n^*) \tilde\phi(U)| \\
&<& |\tilde\phi (B_n^* \gamma^L(U) )-  \tilde\phi(B_n^*) \tilde\phi(U)|
 + \varepsilon\\
&=& | \langle \gamma^{-L}(B_n) \xi , U\xi \rangle -
\tilde\phi(B_n^*) \tilde\phi(U) | + \varepsilon \\
&\le& |\langle \gamma^{-L}(B_n) \xi , U\xi \rangle - 
\langle \gamma^{-L}(B_n) \xi , U_m \xi \rangle |  \\
&&+ |\langle \gamma^{-L}(B_n) \xi , U_m\xi \rangle -
\tilde\phi(B_n^*) \tilde\phi(U) |
+  \varepsilon \\
&<& |\tilde\phi (B_n^* \gamma^L(U_m) )-\tilde\phi(B_n^*) \tilde\phi(U)
 | + 2\varepsilon \\
&<& 
|\tilde\phi (B_n^* \gamma^L(U_m) ) -
 \tilde\phi (B_n^*) \tilde\phi(U_m) | 
 +3 \varepsilon \\
&<&  4\varepsilon.
\end{eqnarray*}
Therefore, we have 
\[
\langle B_n \xi , U \xi \rangle = \tilde\phi(B_n^*) \tilde\phi(U).
\]
This means $U \xi = \tilde\phi (U) \xi$.
Since $\tilde\phi$ is faithful, $\xi$ is a separating vector.
Therefore, we have $U = \tilde\phi(U) I$.
This is a contraction if $g$ is not a unit.
\end{proof}

Let $\{ e_{ij} \}$ be a system of matrix units of $M_d$ and
we write $(i_1 i_2 \ldots i_n, j_1j_2 \ldots j_n)$ for
$e_{i_1 j_1} \otimes e_{i_2j_2} \otimes \cdots \otimes e_{i_nj_n}$.

\begin{example}\label{exam1}{\rm
Let $d=3$ and $\psi$ be a state on $M_3$ whose density matrix is given by
\[
D = {\rm diag}({1 \over 2+ \lambda} ,{1  \over 2 +\lambda }
 , {\lambda  \over 2+\lambda })
 ={1 \over 2+\lambda } {\rm diag}(e^{\log 1},
 e^{\log 1 }, e^{\log \lambda })
\]
for $\lambda \ge 0$ with $\lambda \neq 1$.
We take mutually prime numbers $n,m\in{\mathbb N}$ and
define the gauge group $G$ by
\[
G = \{ {\rm diag}(1, e^{2\pi i k \over m} , e^{2\pi i l \over n} )\,\, |
\,\, 1 \le k \le m, 1 \le l \le n \}.
\]
Let $\tilde\phi = \bigotimes_{n= -\infty}^\infty \psi$, then
$\tilde\phi$ is a $C^*$-finitely correlated state and
$\phi = \tilde\phi | {\mathfrak B}^G$ is a generalized Markov chain.
It follows from Proposition \ref{main} that $\pi({\mathfrak B})''$ is
a type ${\rm III}_\lambda$ factor.
By [\ref{ohno2}] and a simple caluculation, we see that
${\mathfrak A}_{[1,n]} = {\mathfrak B}^G$ 
is the linear span of $(i_1\ldots i_n , j_1 \ldots j_n)$
such that $|\{ k \,| \, i_k = 2\}| \equiv |\{l \,|\, j_l =2\} |$
in mod $m$ and $|\{ k \,| \, i_k = 3\}| \equiv |\{l \,|\, j_l =3\} |$
in mod $n$.
For such $(i_1\ldots i_n , j_1 \ldots j_n)$, we have
\[
D^{it} (i_1\ldots i_n , j_1 \ldots j_n) D^{-it}
= e^{it nk \log \lambda} (i_1\ldots i_n , j_1 \ldots j_n)
\]
for some $k \in {\mathbb  Z}$.
Hence, $\pi({\mathfrak A})''$ is a subfactor of 
$\pi({\mathfrak B})''$ with index $nm$ and
of type ${\rm III}_{\lambda^n}$.
In fact, the subfactors of type 
${\rm III}_{\lambda}$ factors were classified in [\ref{loi}].
}
\end{example}

\begin{example}
{\rm
Let $d=2$ and $\phi$ be a Markov state generated by
$(E, \rho)$.
We assume that the range of $E$ is the diagonals of $M_2$ and
$E((11,11)) = \lambda_1 e_{11}$, $E((12,12)) = (1-\lambda_1)e_{11}$,
$E((21,21)) = (1- \lambda_2) e_{22}$ and $E((22,22)) = \lambda_2 e_{22}$.
By Proposition \ref{type1}, 
${\rm sp}(\sigma^\phi)$ is a closed subgroup of ${\mathbb R}$
generated by $\log \lambda_1 - \log \lambda_2$ and
$\log \lambda_1 + \log \lambda_2 - \log (1-\lambda_1) 
- \log (1-\lambda_2)$.

Let 
\[
G = {\mathbb Z}_m =\{ {\rm diag} (1, e^{2 \pi k \over m}) \,\, | \,\,
0 \le k \le m-1 \}.
\]
Then, $E$ satisfies the $G$-covariant condition.
Moreover, by the similar calculation as in Example \ref{exam1}
and Proposition \ref{type1}, we obtain that
${\rm sp}(\sigma^\phi | {\mathfrak A})$ is a 
closed subgroup of ${\mathbb R}$ generated by
$m (\log \lambda_1 - \log \lambda_2)$ and
$\log \lambda_1 + \log \lambda_2 - \log (1-\lambda_1) 
- \log (1-\lambda_2)$.

Let $0<\lambda <1$.
Let $n$ and $m$ be mutually prime numbers and 
$G = {\mathbb Z}_{nm}$.
Then, by choosing $\lambda_1, \lambda_2$ satisfying
$\log \lambda_1 - \log \lambda_2 = \log \lambda$
and $\log \lambda_1 + \log \lambda_2 - \log (1-\lambda_1) 
- \log (1-\lambda_2) = n \log \lambda$,
we can construct
a subfactor $\pi({\mathfrak A})''$ of 
a type ${\rm III}_\lambda$ factor $\pi({\mathfrak B})''$ 
with index $nm$ and
of type ${\rm III}_{\lambda^n}$.
}
\end{example}

%%%%%%%%%%%%%%%%%%%%%%%%%%5%%%%%%%%%%%%%%%%%%%%%%%%%%%%%%%reference

\end{document}